\documentclass[aps,showpacs,
]{revtex4}
\usepackage{graphicx}
\begin{document}

\title{\large \bf
Global Stabilization of Controlled Nonlinear System
 "Inverted Pendulum on a Cart" Using Method of Two Lyapunov Functions}
\author{B. L. Mazov}
\affiliation{N.Novgorod State Technical University, Nizhny Novgorod, 24 Minin str.
603600 Russia}
\begin{abstract}
In this work a mathematical  apparatus for method of several global Lyapunov
functions
is applied to study the stability properties of nonlinear model of concrete
singularly perturbed mechanical system "inverted pendulum on a cart" with
discontinuous relay-type control. In this analysis, the concept of solution in
sense of Gelig et al.\cite{GLY} is used.
\end{abstract}
\maketitle

\section{Introduction}

The nonlinear controlled model system "inverted pendulum on a cart" is one of
intensively studied last time (see, e.g. \cite{KO}). The results of
investigation of this system are applicable at analysis of a number of concrete
problems of stability, e.g. monorail \cite {MNF}, fault tolerant control
\cite{CD}, inverted pendulum on a rotating arm as well as satellites and
underwater vehicle with internal rotors etc.\cite{BlLM}.

In the work \cite{BLS} it was performed a study of stability of nonlinear
system "inverted pendulum on a controlled cart" \, for a number of cases. There
the conrol was continuous and was taken as sum of the controls $$ u =
\sum\limits_i{\lambda_i u_i}, \eqno(1) $$ and it was applied to the cart moving
on the plane surface.

In recent work of Bloch et al. \cite{BlLM}, it was performed the investigation
of stabilization of this mechanical system using the method of controlled
Lagrangian. The control law was obtained in the form $$ u = \frac{k b sin \beta
(a {\dot \beta}^2 + c cos \beta)} {a - b^2/{\bar m}(1+ k) {cos \beta}^2},
\eqno(2) $$ where $b = ml, c = - mgl$. $m$ is the pendulum mass, $l$ is the
pendulum lehgth, $g$ is the acceleration due to gravity, $\bar m = M + m$, $M$
is the cart mass. And equilibrium $\beta = \dot \beta = r = \dot r = 0$ is
achieved if dimensionless constant $$ k > \frac{a{\bar m} - b^2}{b^2} =
\frac{M}{m} > 0 \eqno (3) $$ The control force acts to a cart while there are
no direct action to a pendulum. At analysis it is used a character of the Lie
symmetry group for nonlinear system, which in the case of plane system
"inverted pendulum on a controlled cart" \, appears to be that of translation.

Also in that work it is studied an asymptotic stabilization of the system "inverted
pendulum on a controlled cart", and control under investigation is breaked into
"conservative"  and "dissipative" piece
$$ u =
u^{cons} + u^{dis},
\eqno(4)
$$
which are analysed separately.

Introduced in the work Lyapunov function for system with control is taken so
that its total time derivative was non-negative elsewhere ($dV/dt \le 0$) and
disappeared at the set $M$, which is determined by equating to zero of dissipative
control component $u^{dis} = 0$.

As follows from simulation results for concrete parameters of nonlinear system
"inverted pendulum on a controlled cart" \, (for given stabilizing action with
added dissipation), obtained using the MATLAB system, the pendulum begins the
motion from almost horizontal position what indicates that attraction region is
large enough, even at given positive (down) initial velocity. Initial position
of the cart is $s(0) = 0$, and initial velocity is $\dot s = - 3$ m/sec. The
cart comes to an equilibrium state in the position which is far enough from the
initial one. The control law was taken with an initial peak to provide the
initial large acceleration. The Lyapunov function is negative in the beginning
and then it strictly increases up to zero at the equilibrium.

In present work, it is performed a detailed analysis of the behavior of the
controlled mechanical system  "inverted pendulum on a cart" \, with
discontinuous relay-type control on the basis of the method of investigation of
the global asymptotic stability of nonlinear dynamical systems with using of
two Lyapunov functions. The distinct feature of such analysis from the case of
continuous right-hand sides \cite{BLS} is that total time derivatives of the
Lyapunov functions appears to be also discontinuous. Because of this to study
the global stability of this system it is used the method of two Lyapunov
functions (see \cite{BM992}) which permits to perform the stability analysis
for discontinuous systems, in particular, for the case of mechanical systems
with dry friction (see, e.g. \cite{AVH}). At such analysis the solution is
considered in the sense of Gelig et al. \cite{GLY}.

\section{Dynamics of control object in new variables}.  \\[3mm]

It is considered the problem of global stabilization of inverted pendulum on a
controlled cart in the presence of unknown disturbance acting at cart. This
mechanical dynamic system is described by nonlinear system of differential
equations in the form (see, e.g. Brusin et al. \cite{BLS}) $$ \left\{
\begin{array}{lcl}
M\ddot r + L\ddot \beta cos\beta + N\dot r - L\dot\beta^2sin\beta = Gu(t)+D(t)
\\
L cos\beta \ddot r+ I \ddot \beta + c\dot\beta + \kappa \dot r cos\beta - Lg
sin \beta=0,\\
\end{array}
\right.
\eqno (5)
$$
where $D(t)$ is uniformly bounded ($ D(t)\le \Pi$ where $\Pi$ is known value),
continious function of $t$, describing the external action at cart;
$\bar m =M + m$, $L=ml$,
$
c=c_0+\kappa l,$ $I=J+m l^2
$
where $M > 0, m > 0$  are the cart and the pendulum masses,
respectively. \, $N > 0$ and $\kappa > 0$ are
the coefficients of friction force for motion of cart and pendulum, respectively;
$c_0$ is the
coefficient of elastic force moment for rotation friction of the pendulum;
$g$ is the acceleration due to gravity;
$G > 0$ is the coefficient of amplification of motor; $J > 0$ is the inertion
moment for pendulum relative to the mass center; $r$ is the coordinate of cart mass
center; $|\beta| < \pi/2$ is the angle between pendulum axis and the vertical,
numbered from vertical unstable positon of pendulum; $u$ is the value of
controlled signal of regulator.

The purpose of control is to reduce controlled "cart" in asymptotics to given
position and attached "pendulum" in vertical position (at which its gravity
center is situated above fixed point) from any initial position (such that
$|\beta| < \pi/2)$ of "cart" and "reversed pendulum" in the presence of
immeasurable perturbation  at cart, i.e. $\beta \to 0, r \to 0 \quad as \quad t
\to \infty.$\\

Let introduce new variables $(s, \dot s, \Omega, \dot \Omega)$ for state vector
of the object of control using relations
$$ s = r + \varrho \Omega \, (\beta),\quad
\Omega (\beta)= -\ln{|\tan(\frac{\pi}{4}-\frac{\beta}{2})|}, \, \beta \in
(-\frac{\pi}{2}, \frac{\pi}{2}), \, \rho > 0.
\eqno(6)
$$
$$
\beta = 2(\pi
/4 - arctg(exp(-\Omega))).
\eqno(7)
$$

The variables $(s, \dot s)$ are connected due to system (5) by the following
equation \cite{BLS} $$ (p + a)[\frac{1}{A}\psi(\beta)(s + \dot s) + bs] =
\sum^{11}_{i = 1} \lambda_i u_i + b(\dot s + \alpha s) - u - \frac{D(t)}{G},
\eqno(8) $$ $$ \psi(\beta) = MI - L^2{cos}^2\beta, $$ where $p = d/dt; \, a >
0; \, \alpha > 0; \, b \ge 0$ are some numbers. $A = |G(I - \rho L)| > 0;
\, u_i$ are the continuous functions of new variables
\\ $$ \left\{
\begin{array}{lcl}
u_1 = \dot r, \, u_2 = \dot r {cos}^2\beta, \, u_3 = {\dot \beta}^2 sin{\beta},
\, u_4 = \dot \beta cos\beta , u_5 = sin 2\beta, \\ u_6 = \dot \beta/cos\beta =
\dot \Omega, \, u_7 = tg\beta, \, u_8 = a\dot s + \rho {\dot
\beta}^2sin\beta/cos^2\beta, \\ u_9 = (\dot \beta sin 2\beta - a
cos^2\beta)\dot s - \rho {\dot \beta}^2 sin\beta, \, u_{10} = \dot s + as, \\
u_{11} = (p + a)(scos^2\beta) = (\dot s + as)cos^2\beta - ssin2\beta.
\end{array}
\right. \eqno(9) $$
\\ Here functions of old variables are expressed via new ones with taking into
account the relations directly following from (6) $$ \dot r = \dot s - \rho
\dot \Omega (\beta), \quad \dot \beta = \dot \Omega cos\beta, \eqno(10) $$ and
angle $\beta$ is expressed from (7).

The constants $\lambda_i$ is in the form  \\
$$
\begin{array}{lcl}
\displaystyle \lambda_1 = \frac{\rho M\kappa + IN - \rho N L}{A}; \quad
\lambda_2 = -\frac{L\kappa}{A}; \quad \lambda_3 = \frac{L}{G}; \\
\displaystyle \lambda_4 = -\frac{Lc}{A}; \quad \lambda_5 = -\frac{L^2g}{2A};
\quad \lambda_6 = -\frac{\rho Mc}{A}; \quad \lambda_7 = -\frac{\rho LMg}{A} \\
\displaystyle \lambda_8 = -\frac{MI}{A}; \quad \lambda_9 = -\frac{L^2}{A};
\quad \lambda_{10} = -\frac{MI\alpha}{A}; \quad \lambda_{11} = \frac{\alpha
L^2}{A} \\
\end{array}
\eqno(11) $$                                \\ Further, if the relations $$
\ddot r = \ddot s -\rho \ddot \Omega, \quad \ddot \beta = \ddot \Omega cos\beta
- {\dot \Omega }^2 sin\beta cos\beta, $$ following from (10), will be taken
into account, then equation for $\Omega$ can be expressed in new variables as
\cite{BLS} $$ \ddot \Omega + d_1\dot \Omega + d_2 tg\beta + d_3 {\dot \Omega}^2
sin\beta = \frac{L\ddot s + \kappa \dot s}{\rho L - I} \eqno(12) $$ where new
designations $$ d_1 = \frac{\kappa \rho - c}{\rho L - I}, \quad d_2 =
\frac{Lg}{\rho L - I}, \quad             \, d_3 = \frac{I}{\rho L - I}
\eqno(13) $$ are introduced. The obtained equations (8) and (12) describe the
behavior of initial system in new variables.

Let take the control law in the form  \cite{BLS}
$$ u(t) = \bar u(t) + \Delta u(t),
\eqno(14)
$$
where function \cite{BLS}
$$ \bar u(t) = \sum^{11}_{i =
1}{\lambda_iu_i(t)}
\eqno(15)
$$
is taken so that right-hand side of (8) was equal to zero
if $\Delta u(t) = 0$. $\Delta u(t)$ - correcting term to compensate perturbation
$D(t)$, which is a piecewise function with a jump discontinuity along some smooth
suface.

Then (8) for closed system at $b = 0$ will be wrote in the form (see
\cite{BLS}) $$ (p + a)[\frac{1}{A}\psi(\beta)(\dot s + \alpha s)] = - \Delta
u(t) - \frac{D(t)}{G}. \eqno(16) $$ Let denote $$ \psi(\beta)(\dot s +
\alpha s) = \gamma(t), \eqno(17) $$ then obtain equation for $\dot s$ $$ \dot
s = \frac{\gamma}{\psi(\beta)} - \alpha s . \eqno(18) $$

Substitution of (17) in (16) results in  $$ (p + a)\frac{\gamma}{A} = -\Delta u(t)
- \frac{D(t)}{G} \eqno(19) $$ или $$ \dot \gamma = - a\gamma - A\Delta u(t)
- B D(t), \eqno(20) $$ где $B = A/G$.

The equations for $y = col(\Omega, \dot \Omega)$ (see (12) can be obtained, if
express $\ddot s$ from (18) $$ \ddot s = \frac{\dot \gamma \psi - \gamma \dot
\psi}{{\psi}^2} - \alpha \dot s \eqno(21) $$ or, with taking into account (18)
$$ \ddot s = \frac{\dot \gamma \psi - \gamma \dot \psi}{{\psi}^2} -
\alpha(\frac{\gamma}{\psi} - \alpha s) \eqno(22) $$ For numerator of right-hand
side of (12) with taking into account (21) one obtains $$ L\ddot s + \kappa
\dot s = \frac{L}{{\psi}^2}(\dot \gamma \psi - \gamma \dot \psi) + (\kappa -
\alpha L)\dot s \eqno(23) $$ or with taking into account (20) $$ L\ddot s +
\kappa \dot s = \frac{L}{{\psi}^2}((-\alpha \gamma - A\Delta u(t) - B D(t))\psi
- \gamma \dot \psi) + {\dot s}(\kappa - \alpha L) \eqno(24) $$ Substituting
(24) in (12) one obtains $$
\begin{array}{lcl}
\ddot \Omega = - d_1\dot \Omega - d_2 tg\beta - d_3 {\dot \Omega}^2 sin\beta +
\\ \displaystyle + \frac{1}{\rho L - I}(L\frac{(-a \gamma - A\Delta u(t) - B
D(t))\psi - \gamma \dot \psi}{{\psi}^2} + {\dot s}(\kappa - aL)
\end{array}
\eqno(25) $$ где $\dot \psi = - L^2 sin 2\beta$.

So, the dynamics of the control object (5) in new variables
$s, \gamma, \Omega, \dot{\Omega}$ will be determined by system of four differential
equitions of the first order  \\ $$ \left\{
\begin{array}{lcl}
\displaystyle \frac{ds}{dt} = \frac {\gamma}{\psi (\beta)}-\alpha s \\
\displaystyle \frac{d\gamma}{dt} =-a \gamma -A \triangle u(t)- B D(t) \\ [3mm]
\displaystyle \frac{d\Omega}{dt} = \dot \Omega \\ [3mm]
\displaystyle \frac{{d\dot{\Omega}}}{dt} = - d_1\dot \Omega - d_2 tg\beta - d_3
{\dot \Omega}^2 sin\beta +  \\ \displaystyle + \frac{L}{\rho L - I}\frac{(-a
\gamma - A\Delta u(t) - B D(t))\psi - \gamma \dot \psi}{{\psi}^2} + {\dot
s}(\kappa - aL) \\
\end{array}
\right. \eqno (26) $$

\section{The case of continuously acting perturbation $D(t)$}.  \\[3mm]

To stabilize the system it will be used a discontinuous relay-type function $$
\triangle u(t)=\bar \Pi sign \,\gamma, \, sign \,\gamma=\\
\left\{\begin{array}{lcl} \displaystyle {}1, \, \gamma>0\\
 -1, \, \gamma<0\\
(-1, 1), \, \gamma=0\\
\end{array}\right.
\eqno (27) $$                      \\ [1mm] So, as a control it will be used a
function $$ u = \bar u + \Delta u(t), $$ where $\bar u$ is taken from (15).

If denote $x=col(s,\gamma),\, y=col(\Omega,\dot \Omega)$, then system
(26) can be attributed to the class of systems considered in \cite{BM992}
.

To use theorems 1-3 from \cite{BM992} there is necessary the verification of all conditions
and proposals used.

It can be obtained that conditions 1.2, 1.3 of theorem 1 from \cite{BM992} are fulfilled
if functions  $V,W,\psi,\eta$ are taken in the form \\ [1mm] $$ \left\{
\begin{array}{lcl}
V(s,\gamma )=s^2+k \gamma^2,\, \\ \displaystyle W(\Omega, \dot \Omega) =\frac
{\dot\Omega^2}{\cos^{q-2} \beta(\Omega)}+ \displaystyle \frac r{\cos^{q-1}
\beta(\Omega)}-r\, > \, 0 \\ \psi =\varepsilon (|s|^2+|\gamma|^2) \\
\displaystyle \eta=\frac{2a \dot\Omega^2}{\cos^{q-2}\beta(\Omega)}
\end{array}
\right. \eqno (28) $$              \\ [1mm] where
$
k>0, \,\varepsilon >0
$
are the numbers taken by corresponding manner \,
$
q=2\varrho L(\varrho L-I)^{-1}>1, \, r=2(q-1)^{-1}Lg(\varrho L-I)^{-1}>0. $ \\

\section{The case of immeasurable velocities: reducing to
singularly perturbed system}

Consider system (5) at $D(t) \equiv 0$. Let that in the regulator
velocities $\dot r, \dot \beta$ are measured inexactly, and their approximate
values $\hat z_1, \hat z_2$ are generated with using the system
\cite{Bru951,Bru952}. \\ [1mm] $$ \left\{
\begin{array}{lcl}
\displaystyle \mu \frac{{d\hat z_1}}{dt} = \dot r - \hat z_1,
\\[3mm] \displaystyle \mu \frac{{d\hat z_2}}{dt} = \dot \beta -
\hat z_2 \\
\end{array}
\right. \eqno (29) $$                                   \\ [1mm] where $0 < \mu
< 1$ is the small parameter. The processes $\hat z_{1, 2}$ approximate
processes $\dot r, \dot \beta$. It is essentially that in contrast to later
they can be calculated without using of differentiation operation $$ \hat
z_1(t) = D_{\mu}^1(r(t)), \quad \hat z_2(t) = D_{\mu}^1(\beta(t)). \eqno(30) $$
In this case instead of globally stabilizing control
 $u = \sum\limits_i^{}\lambda_i u_i(t)$
consider control in the form \\ [1mm]
$$ u =
\sum\limits_i^{}{\lambda_i \hat {u_i}(t)}
\eqno(31)
$$
\\ [1mm]
where $\hat {u_i}$ is determined from continuous functions
$u_i$ by replacement $\dot r, \dot \beta$ to $ z_1, z_2$, respectively. $$
\left\{
\begin{array}{lcl}
\displaystyle u_1 = \hat z_1, \, u_2 = \hat z_1 {cos}^2\beta, \, u_3 = {\hat
z_2}^2 sin\beta, \, u_4 = \hat z_2 cos\beta , \\ \displaystyle u_5 = u_5, \,
u_6 = \hat z_2/cos\beta = \dot \Omega, \, u_7 = u_7, \, u_8 = a\hat z_3 +
\rho{\hat z_2}^2sin\beta/cos^2\beta \\ \displaystyle u_9 = (\hat z_2 sin 2\beta
- a cos^2\beta)\hat z_3 - \rho {\hat z_2}^2 sin\beta; u_{10} = \hat z_3 +
\alpha r - \alpha \rho ln|tg(\frac {\pi}{2} - \frac{\beta}{2}|; \\
\displaystyle u_{11} = \hat z_3 cos^2\beta + s\hat z_2sin2\beta + ar - a\rho
ln|tg(\frac{\pi}{2} - \frac{\beta}{2})|cos^2\beta.
\end{array}
\right. \eqno(32) $$ Let introduce a vector process $z(t) = col(z_1(t),
z_2(t))$ $$ z_1 = \hat z_1 - \dot r, \quad z_2 = \hat z_2 - \dot \beta
\eqno(33) $$ and also denote $$ u^{\mu}_i = \hat u_i - u_i . \eqno(34) $$ Then
system (26), (29), (31) with taking into account (33), (34) can be presented in
the form of six diferential equitions of the first order (cf. with (5)) \\ $$
\left\{
\begin{array}{lcl}
\displaystyle \frac{ds}{dt} + \alpha s = \gamma / \psi(\beta ) \\ [3mm]
\displaystyle \frac{d\gamma}{dt} + a\gamma = -\sum\limits_{i=1}^n \lambda_i
u_i^{\mu} \\ [3mm] \displaystyle \frac{d\Omega}{dt} = \dot \Omega \\  [3mm]
\displaystyle \frac{d\dot \Omega}{dt} = f_\Omega (s, \gamma, \mu), \\ [3mm]
\displaystyle \mu \frac{d{z_1}}{dt} = - z_1 - \mu f_r (s,\gamma , \Omega, \dot
\Omega,z_1,z_2) \\        [3mm] \displaystyle \mu \frac{d{z_2}}{dt} = - z_2  -
\mu f_\beta (s,\gamma , \Omega, \dot \Omega,z_1,z_2) .
\end{array}
\right. \eqno(35) $$                                 \\ [1mm]
From above presented it's seen that
 $f_\Omega, f_r, f_{\beta}$ are the continuous functions
. Moreover, at $\mu=0$ system (35) is reduced to system (26)
(at $D(t)=0$ and $\triangle u = 0$).

To obtain functions $f_r$ and $f_{\beta}$ in direct form express $\hat z_1$,
$\dot {\hat z_1}$ и $\hat z_2$, $\dot {\hat z_2}$ via $z_1$ и $z_2$ $$
\begin{array}{lcl}
\hat z_1 = z_1 + \dot r, \quad \dot {\hat z_1} = \dot z_1 + \ddot r, \quad \hat
z_2 = z_2 + \dot{\beta}, \quad \dot {\hat z_2} = \dot z_2 + \ddot \beta .
\end{array}
\eqno(36) $$ Substituting these expressions in (4.5.1) one obtains $$ \mu (\dot
z_1 + \ddot r) = - z_1, \quad \mu (\dot z_2 + \ddot \beta) = - z_2 \eqno(37) $$
from where $$ \left\{
\begin{array}{lcl}
\displaystyle \mu \frac{dz_1}{dt} = - z_1 - \mu \ddot r \\ [3mm]
\displaystyle \mu \frac{dz_2}{dt} = - z_2 - \mu \ddot \beta \\
\end{array}
\right. \eqno(38) $$ So, for vector process (33) $$ f_r (s,\gamma , \Omega,
\dot \Omega,z_1,z_2) = \ddot r, \quad f_\beta (s,\gamma , \Omega, \dot
\Omega,z_1,z_2) = \ddot \beta \eqno(39) $$ Express values $\ddot r$ and $\ddot
\beta$ (due to (5)) as continuous functions of state of the control object and
regulator \cite{Bru952}. Let write the system (5) in the form $$ \left\{
\begin{array}{lcl}
\displaystyle M\ddot r + L\ddot \beta cos\beta = Gu(t)+D(t) - N\dot r +
L\dot\beta^2sin\beta \\ \displaystyle L cos\beta \ddot r+ I \ddot \beta = Lg
sin \beta - c\dot\beta - \kappa\dot r cos\beta .\\
\end{array}
\right. \eqno (40) $$ The solution of this system of linear algebraic equitions
relative to variables $\ddot r$ and $\ddot \beta$ has a form $$ \left\{
\begin{array}{lcl}
\displaystyle \ddot r = \frac {I(Gu(t) -B\dot r + Lsin\beta{\dot \beta}^2 +
D(t)) - Lcos\beta(Lgsin\beta {\dot \beta} - c{\dot \beta} - b\dot r cos\beta)}
{MI - L^2{cos\beta}^2} \\
\displaystyle \ddot \beta = \frac {M(Lgsin\beta -c{\dot \beta} - b\dot r
cos\beta) - Lcos\beta(Gu(t) - B\dot r + Lsin\beta{\dot \beta} + D(t))} {MI -
L^2{cos\beta}^2}
\end{array}
\right. \eqno(41) $$ In new variables the values $\ddot r$  and $\ddot \beta$
will be expressed with taking into account (10) as $$ \ddot r = \ddot s - \rho
\ddot \Omega, \quad \ddot \beta = \ddot \Omega cos\beta -{\dot \Omega}^2
sin\beta cos\beta. \eqno(42) $$ где $$ \ddot s = \frac{\dot \gamma \psi -
\gamma \dot \psi}{{\psi}^2} - \alpha \dot s \eqno(43) $$ In correspondence with
(26) $$ \dot \gamma = - a\gamma - A\triangle u(t) - BD(t) \eqno(44) $$ from
where $$
\begin{array}{lcl}
\displaystyle \ddot \Omega = d_1\dot \Omega - d_2 tg\beta - d_3 {\dot \Omega}^2
sin\beta + L\frac{-\alpha \gamma - A\triangle u(t) - BD(t) + 2 \gamma L^2 sin
\beta}{(\rho L - I){\psi}^2} - \\
\displaystyle (\frac{\gamma}{\psi(\beta)} + \alpha s)( aL)
\end{array}
\eqno(45)
$$
Here $\psi$ and its derivative are determined with using of relations
$$
\psi = MI - L^2 {cos\beta}^2, \quad \dot \psi = 2L^2 cos\beta
sin\beta = L^2sin2\beta
\eqno(46)
$$
and $\beta$, in its turn is determined by relation (7).

Now, turn out to verification of fulfilment of conditions of theorem 3 from
\cite{Bru951}. To present the system considered in the form of \cite{Bru951},
replace (38), (39) to system $$ \left\{
\begin{array}{lcl}
\displaystyle \frac{dz_1}{dt} = - z_1 - \mu_2 f_r \\ [3mm] \displaystyle
\mu_1\frac{dz_2}{dt} = - z_2 - \mu_2 f_{\beta} \\
\end{array}
\right.
\eqno(47)
$$
So, system (35), (38) with functions
$$ \left\{
\begin{array}{lcl}
f_r = \ddot r = \ddot s(s, \gamma, \Omega, \dot \Omega) - \rho \ddot \Omega(s,
\gamma, \Omega, \dot \Omega) , \\ f_{\beta} = \ddot \beta = \ddot \Omega(s,
\gamma, \Omega, \dot \Omega)cos\beta - {\dot \Omega}^2 sin\beta cos\beta
\end{array}
\right.
\eqno(48)
$$
can be reduced to the form \cite{Bru951}, if denote  $x
= col(s, \gamma, \Omega, \dot \Omega)$ и $z = col(z_1, z_2)$. $$ x =
\left(\begin{array}{lcl} \dot s \\ \dot \gamma \\ \dot \Omega \\ \ddot \Omega
\end{array}\right)
\quad \bar f(x, z) = \left(\begin{array}{lcl} \gamma /\psi(\beta) - \alpha s \\
-a\gamma - \sum\limits_{i = 1}^{11}\lambda_iu_i^{\mu} \\ \tilde \Omega = \dot
\Omega \\ \ddot \Omega \,  from  \, (45)
\end{array}\right)
\quad \bar g(x, z) = 0 \eqno(49) $$ $$ \dot z = \left(\begin{array}{lcl}
\dot z_1 \\ \dot z_2,
\end{array}\right)
\quad - \Gamma z = \left(
\begin{array}{lcl}
- z_1 \\ - z_2
\end{array}
\right) \quad \bar h(x, z) = \left(
\begin{array}{lcl}
\ddot s - \rho \ddot \Omega \\ \ddot \Omega cos\beta - \dot \Omega sin\beta
cos\beta
\end{array}
\right)
\eqno(50)
$$

At $\mu_2 = 0, z = 0$ it obtained the initial unperturbed system studied in
\cite{Bru93}. As it was demonstrated in that work this system will have
the dissipativity domain in the form of $D = \{(s, \gamma ,
\Omega , \dot\Omega), V(s, \gamma , \Omega , \dot\Omega) \le C\}$, where $V(s,
\gamma , \Omega , \dot\Omega)$ is the continuous piecewise positive definite
in ${\bf R}^m$ (in sense \cite{Rouch}) function, satisfacting the condition
a) of theorem 3 in  \cite{Bru951}), i.e.
$$
\left.\frac{dV(s,
\gamma , \Omega , \dot\Omega )} {dt}\right\vert_{(s,\gamma) \in S(R)} \le -
\varepsilon , \quad V(s, \gamma , \Omega , \dot\Omega) \to \infty \quad at
\quad |s, \gamma , \Omega , \dot\Omega| \to \infty, $$ where $S(R)$ is the
bounded closed (hyper)surface: $$ \{(s, \gamma , \Omega ,
\dot\Omega), \, V(s, \gamma , \Omega , \dot\Omega) = R\}, $$ and $R, \varepsilon$
are the positive numbers.

Further, for system (38) $$ \left\{
\begin{array}{lcl}
\displaystyle \mu \frac{dz_1}{dt} = - z_1 - \mu (\ddot s - \rho \ddot \Omega),
\\ [3mm] \displaystyle \mu\frac{dz_2}{dt} = - z_2 - \mu (\ddot \Omega cos\beta
- {\dot \Omega}^2 sin\beta cos\beta) \\
\end{array}
\right. \eqno(51) $$ where $\ddot s$  and $\ddot \Omega$ are determined from
(21), (25), it will take place a condition  b) of the theorem 3 from
\cite{Bru951}), and as $W(z_1, z_2)$ \cite{Bru951} it can be taken quadratic
form $$ W = <z, Bz>, \quad \Gamma^TB + B\Gamma = - E \eqno(52) $$ i.e. $$
\frac{dW(z_1, z_2)}{dt} \le - \frac{L_1}{\mu_1}|z|^2, \quad (\exists
\bar\mu_1), \, (\forall \mu_1 \in (0, \bar\mu_1]) \quad (W(z) \to \infty \quad
при \quad |z| \to \infty). \eqno(53) $$ Also from  (32), (33) it will follow a
fulfilment of condition c) of theorem 3 in \cite{Bru951}. I.e. all conditions
of theorem 3 will be hold and thus the system (5), (31)-(33), (38) will have
the dissipativity domain $D = \{y, W_{\rho}(x, z) \le C\}$, independing of $\mu
\in (0, \bar \mu]$  \, $\exists \bar\mu > 0$, где $W_{\rho}(y) = V(s, \gamma ,
\Omega , \dot\Omega) + \rho W(z_1, z_2)$ is the global Lyapunov function
satisfacting due to the system considered the equality \cite{Bru951} $$
\left.\frac{dW_{\rho}(y)}{dt}\right\vert_{\mu_2 = 0} = - (|x|^2 + \rho |z|^2).
$$ Here $y = (x, z), \, |y|^2 = |(x, z)|^2 = |x|^2 + |z|^2$, $C$ и $\rho$ are
some positive numbers.

Then, it can be stated the following theorem.

{\bf Theorem }.{\it Given system has attraction domain of equilibrium state
$x=0,z=0$ in the form $W(x,z) \le C$ at all $0 < \mu \le \bar \mu_1$,
and $ C \to \infty $ при  $\bar \mu_1 \to 0$}.

The proof is based on estimation of total time derivative of above introduced
function $$ W_{\rho}(x, z) = W_{\rho}(s, \gamma , \Omega , \dot\Omega ; z_1,
z_2) = V_{\rho}(s, \gamma , \Omega , \dot\Omega) + \rho W(z_1, z_2) $$ with
using of inequality $W_{\rho}(x, z) \le C$ as well as following from it
inequality $|(s, \gamma , \Omega , \dot\Omega)| \le R_{W_{\rho}}(C)$, where
$R_{W_{\rho}}(C) = sup|(s, \gamma , \Omega , \dot\Omega)|, \, W_{\rho}(s,
\gamma , \Omega , \dot\Omega; z_1, z_2) \le C$ following from conditions of the
theorems estimations of right-hand sides $$ \bar f(s, \gamma , \Omega , \dot
\Omega; z_1, z_2) \le
\\
\displaystyle const(C) |(s, \gamma , \Omega , \dot \Omega)| |(z_1, z_2)| $$, $$
\bar h(s, \gamma , \Omega , \dot \Omega; z_1, z_2) \le \\ \displaystyle
const(C)|(s, \gamma , \Omega , \dot \Omega)||(z_1, z_2)| $$ correct in this region.

\section{Results of numerical simulation}

In work \cite{BLS} it was performed a numerical simulation for the case of
closed-loop system with continuous control, described by equation (5). The
behavior of mass-center cart coordinate as well as angular coordinate with
increase of $t$ evident about stabilization of the system for a short enough
time interval at any values of parameters used. Moreover, with increase of
parameter $a$ oscillatory regime of decrease of  angular as well as coordinate
amplitude is changed to damping one.

In present work, the numerical simulation of the behavior of the nonlinear
system "inverted pendulum on a cart" with relay-type control under acting of
continuously acting immeasurable perturbation was performed on the basis of
MATLAB system. In particular, it was performed a numerical simulation of
nonlinear dynamical system on the basis of the system (26) where as control
force it was taken discontinuous relay-type (27) function with following
parameter values   $\alpha = 1; \, a = 0.5; \, k = 0.1; \, I = 1; \, L = 1; \,
M = 1.5; \, d_1 = d_2 = d_3 = 1; \, \rho_m = 2$. The initial values of
variables were the following $s = - 0.7; \, \gamma = 0.7; \, \Omega = 1.0; \,
\dot \Omega = 0.5$.

For the case of absence of the control and external action ($A = 0, \, B = 0$)
in the system all time dependencies for four variables $s, \gamma , \Omega ,
\dot \Omega $ are damping and asymptotically approaching zero but the character
of approaching zero appears to be essentially different. The value $s$ coming
from the point with negative coordinate then becomes to be positive and
approach zero from the positive value side. In contrast, the variable $\gamma$
is damping with increasing time up to zero being positive value beginning from
$t = 0$ ($\gamma (0) = 0.7$). Then, time dependences $\Omega(t)$ and $\dot
\Omega (t)$ are of oscillatory damping character. At $t > 10$ the system comes
to equilibrium state $s =0, \, \gamma = 0, \, \Omega = 0, \, \dot \Omega  = 0$
and rate of approaching the equilibrium state is determined by the friction
force given by parameters $L$ and  $\kappa$ in the nonlinear system (35).

The appearance of discontinuous relay-type control (27) in the system ($A =
0.03$) in absence of external action ($B = 0$) essentially decreases time of
approach the equilibrium state in the system, moreover in $\gamma (t)$
dependence there appear horizontal regions $\gamma = 0$ corresponding to the
sliding regime at the breaking surface in the state space. Note also that
presence of relay-type control leads to significant decrease of time necessary
for coming of the system in the equilibrium state $\beta = 0, \, r = 0$.

The numerical simulation in the case of external immeasurable sinusoidal
perturbation in the system $D(t) = sin t$ in absence of the control ($A = 0$)
all time dependences don't approach the stationary solution of the system (4.3.21)
$x = 0$ but are of oscillatory character with amplitude ($\Delta x = 0.1$) and
period equal to that of external force $D(t)$. The presence of discontinuous
relay-type control (27) ($A = 0.05$) permit approach the stationary solution
with sufficient accuracy ($\Delta x = 0.03$) already at $t > 10$. Such stabilization
is determined, in essential degree, by appearance at $\gamma (t)$ dependence of
horizontal regions corresponding to sliding regime at the breaking surface
$\gamma = 0$ in the state space of the system (26) what don't permit to variables
reach sufficient oscillation amplitude for time interval between horizontal regions.

The presence of discontinuous relay-type control (27) in the case of action of
external sinusoidal perturbation leads to a fast stabilization of the system
guiding the trajectory at the beaking surface $\gamma = 0$ of the system (26)
and further it does not permit this trajectory essentially decline from this surface
at periodical sign change in external perturbation $D(t)$.  The essential fact here
is that at more high-frequency perturbation the stabilization of the system occurs
for a shorter time as compared with low-frequency (or even constant) external action.

\section{Conclusion}

In this work it was performed a study of global asymptotic stability of
nonlinear mechanical system "inverted pendulum on a controlled cart" in
conditions of continuously acting perturbation with using method of two
Lyapunov functions. It was demonstrated that at given choice of Lyapunov
functions and comparison functions the conditions of the basic theorem of the
method are fulfilled, i.e. system "inverted pendulum on a controlled cart" \,
with discontinuous relay-type control, in conditions of continually acting
perturbation is globally asymptotically stable. Here the solution is considered
in the sense of Gelig et al. \cite{GLY}. There is established the robustness
property of above obtained stabilization algorithms relatively nonideality of
measurements of velocities of both cart and pendulum, since in practice there
are no strict signal differentiators (instead of this the approximative
velocity measurements are used). Also, using the method of two Lyapunov
functions it is demonstrated that at given algorithm but with replace of ideal
derivatives, the stability region persists, and its size at $\mu \to 0$
increases unlimitedly and in the limit coincides with the full phase space.



\begin{thebibliography}{99}
\bibitem{KO}
Kolesnikov A.A. Synergetic Theory of Control // Moscow, Energoatom Press, 1994
\bibitem{MNF}
Mori S., Nishihaca H., Furuta K., Control of unstable mechanical system control
of pendulum Int.J.Control. v.23, 5, 1976
\bibitem{CD}
Cardoso A., Dourado A., Robust model-based tolerant control of a mobile
structure - application to an inverted pendulum // in: Proc.of the 7th
Int.Symp.on Intelligent Robotic Syst., July 20-23, 1999, Coimbra, Portugal.
\bibitem{BlLM}
Bloch A.M., Leonard N.E., Mardsen J.E., Controlled Lagrangians and the
stabilization of mechanical systems I: the first matching theorem // IEEE
Trans.Automat.Control v.45, No.12, p.2253-2270, 2000.
\bibitem{BLS}
Brusin V.A., Leibo A.M., Serebryakov D.K. Global stabilization of instable
nonlinear two-mass system // Izv.of RAS: Tech.Cibernetics No.4, p.3-12, 1991.
\bibitem{AVH}
Andronov A.A., Vitt A.A., Khaikin C.E. Theory of oscillations // Moscow,
Fizmatgiz, 1959 - 916 p.
\bibitem{BM992}
Brusin V.A., Mazov B.L. The method of two Lyapunov functions in the global
stabilization problem for nonlinear systems // Differential Equations
V.35, No.5, p.626-633, 1999
\bibitem{GLY}
Gelig A.Kh., Leonov G.A., Yakubovich V.A. Stability of Nonlinear Systems
with Non-Single Equilibrium State // Moscow, Nauka Press, 1978. - 400 p.
\bibitem{Bru951}
Brusin V.A. On the one class of singularly perturbed adaptive systems I //
Authomatics and Telemechanics No.4, p.119-129, 1995.
\bibitem{Bru952}
Brusin V.A. On the one class of singularly perturbed adaptive systems II //
Authomatics and Telemechanics No.5, p.103-1113, 1995.
\bibitem{Bru93}
Brusin V.A. Global stabilization of the system "inverted pendulum on a cart"//
Izv.of RAS: Tech.Cibernetics No.4, p.30-39, 1993
\bibitem{Rouch}
Rouche N., Habets P., Laloy M. Stability Theory by Lyapunov's Direct Method
Springer-Verlag New York, Heidelberg, Berlin, 1977
\end{thebibliography}
\end{document}